\documentclass[12pt,reqno]{amsart}
\usepackage{fullpage}
\begin{document}
\thispagestyle{empty}

\newcommand{\C}{\mathcal C}
\newcommand{\F}{\mathcal F}
\newcommand{\G}{\mathcal G}
\newcommand{\Q}{{\mathbb Q}}
\newcommand{\R}{{\mathbb R}}
\newcommand{\Z}{{\mathbb Z}}
\newcommand{\floor}[1]{\lfloor#1\rfloor}
\newcommand{\op}[2]{\langle#1,#2\rangle}

\newtheorem{theorem}{Theorem}
\newtheorem{lemma}[theorem]{Lemma}
\newtheorem{prop}[theorem]{Proposition}

\title{The Unreasonable Effectualness of Continued~Function~Expansions}
\author{Greg Martin}
\address{Department of Mathematics\\University of British
Columbia\\Room 121, 1984 Mathematics Road\\Vancouver, BC V6T 1Z2}
\email{gerg@math.ubc.ca}
\maketitle


\section{Introduction}

The familiar continued fraction expansion of a real number has great
importance in its approximation by rational numbers, and the predictable
behavior of the continued fractions of certain classes of real numbers has
added benefits. For example, the fact that the continued fraction expansion of
a rational number terminates is essentially a reexpression of the Euclidean
algorithm; also, the periodicity of the continued fractions for quadratic
irrationals is crucial for calculating the fundamental units of real quadratic
fields. Already in 1848 Hermite, in correspondence with Jacobi, asked about the
existence of generalizations of continued fractions such that algebraic numbers
of given degree would have periodic expansions. Since that time, myriad
different generalizations have been studied (see \cite{schweiger} for an
extensive list). Herein we focus on the $f$-expansions introduced by
Bissinger \cite{bissinger}, which we define momentarily. The purpose of this
paper is to demonstrate that the function $f$ can be chosen so that the
expansions of prescribed real numbers can have essentially any desired
behavior. The following results, listed in roughly increasing order of
unlikeliness, are representative of what we can prove.

\begin{theorem}
For any two real numbers $x,y\in(0,1)$, there exists a function $f$ such that
the $f$-expansion of $x$ is the same as the usual continued fraction expansion
of~$y$.
\label{oneone.thm}
\end{theorem}

\begin{theorem}
There exists a function $f$ such that the $f$-expansion of any rational or
quadratic irrational terminates.
\label{frage.thm}
\end{theorem}

\begin{theorem}
There exists a function $f$ such that the $f$-expansion of a real number $x$
is periodic if and only if $x$ is a cubic irrational number.
\label{cubic.thm}
\end{theorem}

\begin{theorem}
There exists a function $f$ such that, simultaneously for every integer
$d\ge1$, a real number $x$ is algebraic of degree $d$ if and only if the
$f$-expansion of $x$ terminates with the integer $d+1$.
\label{simultaneous.thm}
\end{theorem}

We remark that Vorono\u\i's algorithm \cite{voronoi} for calculating units in
number fields of degree higher than two is not directly relevant to Theorem
\ref{cubic.thm}, since we are interested in generalizations that
give well-defined expansions for every real number. We also remark that in all
four theorems, the cardinality of the set of functions $f$ satisfying the given
property is that of the continuum, which is the cardinality of the set of all
continuous functions on the real numbers. Finally, we mention an even more
surprising generalization of Theorem \ref{simultaneous.thm}, which we discuss
in more detail later in the paper: there exists a function $f$ such that the
$f$-expansion of every algebraic real number $x$ of degree $d$ terminates with
the integer $d+1$, and the $d+1$ integers directly preceding this final $d+1$
encode the minimal polynomial of~$x$.

Let us describe more precisely the class of expansions we shall consider. The
output of any such expansion will be a sequence in the set $\C = \C_i \cup
\C_t$, where
$$
\C_i = \{ [a_0; a_1, a_2, \dots] \colon \text{each }a_j\in\Z,\, a_j\ge1
\text{ for all } j\ge1 \}
$$
and
$$
\C_t = \{ [a_0; a_1, \dots, a_n] \colon n\ge0,\, \text{each }a_j\in\Z,\,
a_j\ge1 \text{ for all } 1\le j\le n,\, a_n\ge2 \text{ if } n\ge1 \}.
$$
We emphasize that the elements of these sets are formal sequences of integers,
not real numbers; the sets $\C_i$ and $\C_t$ are the infinite and terminating
sequences, respectively. Let $\floor x$ and $\{x\}=x-\floor x$ denote the
greatest integer function and fractional part, respectively, of $x$. Define
$\F$ to be the set of decreasing homeomorphisms from $(1,\infty)$ to $(0,1)$,
that is, the set of all strictly decreasing continuous functions $f$ defined
on $(1,\infty)$ satisfying $\lim_{x\to1^+}f(x)=1$ and $\lim_{x\to+\infty}f(x)
= 0$. Throughout this paper, $f$ will denote a function from the class $\F$
and $\phi$ will denote the inverse of $f$, so that $\phi$ is a decreasing
homeomorphism from $(0,1)$ to $(1,\infty)$.\looseness=-1

We define the {\it expansion function\/} $E_f:\R\to\C$ (sometimes called the
representation function by other authors) as follows. Given $x_0\in\R$, we set
$a_0=\floor{x_0}$. If $x_0$ is not an integer, we set $x_1=\phi(\{x_0\})$ and
$a_1=\floor{x_1}$. If $x_1$ is not an integer, we set $x_2=\phi(\{x_1\})$ and
$a_2=\floor{x_2}$, and so on. Then the value of $E_f(x_0)$ is
$[a_0;a_1,\dots,a_n]\in\C_t$ or $[a_0;a_1,a_2,\dots]\in\C_i$, depending on
whether one of the $x_j$ is equal to an integer. The prototypical example
uses the reciprocal function $r(x)=1/x$, in which case $E_r(x_0)$ is the usual
continued fraction expansion of $x_0$. In general, we call $E_f(x_0)$ the {\it
continued function expansion\/} of $x_0$, or sometimes the {\it continued
$f$-expansion\/} for a specific function $f$. In the terminology of
\cite{schweiger}, these are $f$-expansions of type~A ($f$-expansions of type~B
are formed from increasing functions $f$ and generalize the usual decimal
expansions of real numbers).

In this paper, when we write simply $[a_0;a_1,a_2,\dots]$ we mean the formal
element of $\C$. If we want to refer to the real number whose usual continued
fraction expansion is $[a_0;a_1,a_2,\dots]$, we use the notation
$[a_0;a_1,a_2,\dots]_r$. In general, for any function $f\in\F$ we define an
{\it evaluation function\/} $V_f$ on $C_t$ recursively by setting $V_f([a_0]) =
a_0$ and $V_f([a_0;a_1,\dots,a_n]) = a_0 + f(V_f([a_1;a_2,\dots,a_n]))$. Thus
\begin{equation}
V_f([a_0;a_1,\dots,a_n]) = a_0 + f(a_1 + f(a_2 + \dots + f(a_n) \cdots )),
\label{Vffinitedef}
\end{equation}
which is the continued $f$-expansion of a certain real number. We extend the
definition of $V_f$ to as much of $\C$ as we can by defining
\begin{equation*}
V_f([a_0;a_1,a_2,\dots]) = \lim_{n\to\infty} V_f([a_0;a_1,\dots,a_n])
\end{equation*}
when the limit exists. We shall often write $[a_0;a_1,\dots,a_n]_f$ as a
shorthand for $V_f([a_0;a_1,\dots,a_n])$, thus generalizing the notation
$[a_0;a_1,a_2,\dots]_r$ given above.

Note that $[a_0;a_1,a_2,\dots]_r$ always exists and equals the unique real
number $x$ such that $E_r(x) = [a_0;a_1,a_2,\dots]$, that is, the expansion
function $E_r:\R\to\C$ and the evaluation function $V_r:\C\to\R$ are inverses
of each other. For a general function $f\in\F$, the definitions of $E_f$ and
$V_f$ do imply that the restriction of $V_f$ to $\C_t$ and the restriction of
$E_f$ to $V_f(\C_t)$ are inverses of each other; however, $E_f$ might not be
injective on all of $\R$, or $V_f$ might not be defined on all of $\C$, and so
on.

If the function $f\in\F$ does have the property that $E_f$
and $V_f$ are inverses of each other, we call $f$ a {\it faithful\/} function.
In other words, $f$ is faithful if and only if $E_f$ is bijective and the
limit defining $[a_0;a_1,a_2,\dots]_f$ exists for every element of $C_i$ and
always equals the unique real number $x$ such that $E_f(x) =
[a_0;a_1,a_2,\dots]$. (The list of properties in this last
sentence is probably redundant for characterizing faithful functions, though
we shall not need a more streamlined criterion. Our definition of faithful is
related to what other authors call a {\it valid\/} representation.) In this
terminology, Theorems 1--4 can be stated more precisely using the phrase
``there exists a faithful function $f\in\F$ such that the continued
$f$-expansion $E_f$ of \dots.''\looseness=-1

The idea of our method is to endow $\C$ with a topology that is naturally
related to continued fraction expansions and then to treat the expansion
and evaluation maps $E_f$ and $V_f$ as continuous functions between $\R$ and
$\C$. We describe this topology and begin to explore its consequences in
Section 2. With this foundation, we can make substantial progress by composing
these functions with carefully chosen continuous functions from $\R$ to itself.
This technique, which leads to proofs of Theorems 1--4, is expounded in
Section 3. Finally, in Section 4 we investigate the class of
continued function expansions given by power functions $f(x) = x^{-\alpha}$,
including several numerical examples that partially motivated this paper.

At this point we should confess what the reader might already suspect, that the
functions giving the nice behaviors of Theorems 1--4 are infeasible for actual
computations. Indeed, the existence of such functions is essentially a
consequence of the existence of continuous functions on the interval $(0,1)$
with certain properties. We have chosen the title of this paper, a subtle
variation on the famous phrase ``unreasonable effectiveness'' pioneered by
Wigner~\cite{wigner}, for this reason. Mirriam--Webster's Collegiate Dictionary
contrasts the two words by saying that `effective' in fact ``stresses the
actual production of or the power to produce an effect'', while `effectual'
merely ``suggests the accomplishment of a desired result especially as viewed
after the fact''. We cannot think of a more apt description of these
techniques.\looseness=-1

\section{Topological Preliminaries}

We recall that if $S$ is a set endowed with a total linear ordering, the {\it
order topology\/} on S is defined by declaring the open sets to be arbitrary
unions of open intervals in $S$, that is, of sets of the form $(x,y) =
\{s\in S\colon x<s<y\}$. If $S$ and $T$ are two ordered sets, a function
$f:S\to T$ is {\it order preserving\/} if, whenever $x<y$ in $S$, we have
$f(x)<f(y)$ in $T$. It is easily checked that this implies $x<_S y
\Longleftrightarrow f(x)<_T f(y)$, and consequently any order-preserving
function is automatically injective. We shall need the following result as
well:

\begin{lemma}
Let $h:S\to T$ be a function between the two ordered spaces $S$ and $T$. If
$h$ is order-preserving and surjective, then $h$ is a homeomorphism.
\label{orderbi.lemma}
\end{lemma}

\begin{proof}
Since any order-preserving function is injective, $h$ is in fact a bijection,
and it is easily verified that $h^{-1}$ is also order-preserving. Moreover, it
is true that the image under $h$ of any open interval $(x,y)\subset S$ is
exactly the open interval $(h(x),h(y))$ in $T$: certainly the image is
contained in this open interval by the order-preserving property of $h$, while
every point in $(h(x),h(y))$ must have a preimage in $S$ by the surjectivity of
$h$, and this preimage must be in $(x,y)$ again by order-preservation. This
shows that both $h$ and $h^{-1}$ are continuous, and therefore $h$ is a
homeomorphism.
\end{proof}

We can endow $\C$ with the {\it alternating lexicographic order topology\/},
or the alt-lex topology for short, which is the order topology defined by the
following total ordering of $\C$:
\begin{multline*}
[a_0;a_1,a_2,\dots] < [b_0;b_1,b_2,\dots] \Longleftrightarrow (a_0<b_0)
\text{ or } (a_0=b_0 \text{ and } a_1>b_1) \\ \text{ or } (a_0=b_0 \text{ and }
a_1=b_1 \text{ and } a_2<b_2) \\ \text{ or } (a_0=b_0 \text{ and }
a_1=b_1 \text{ and } a_2=b_2 \text{ and } a_3>b_3) \text{ or } \dots.
\end{multline*}
\goodbreak

\noindent Terminating elements $[a_0;a_1,\dots,a_n]\in\C_t$ are treated as
$[a_0;a_1,\dots,a_n,+\infty]$ when applying this definition. It is easy to see
that for any $f\in\F$, the evaluation map $E_f$ is semi-order-preserving, that
is, $x\le y$ in $\R$ implies that $E_f(x)\le E_f(y)$ in $\C$. In particular, if
$f$ is a faithful function, then the function $E_f$ is bijective and
hence strictly order preserving. We conclude from Lemma \ref{orderbi.lemma}
that the evaluation function $E_f$ of any faithful $f\in\F$ is a homeomorphism
from $\R$ to $\C$.

\begin{lemma}
A subset $B$ of $\C$ is dense if and only if, for every element $x =
[a_0;a_1,\dots,a_n]$ of $\C_t$, there exists an element $b =
[b_0;b_1,b_2,\dots]$ of $B$ such that $b_0=a_0$, $b_1=a_1$, \dots, $b_n=a_n$.
\label{dense.lemma}
\end{lemma}

\begin{proof}
The key observation is that the set of elements $x = [x_0;x_1,x_2,\dots]$ of
$C$ such that $x_0=a_0,\, x_1=a_1,\, \dots,\, x_n=a_n$ is one of the
half-open intervals $\big( [a_0;a_1,\dots,a_n+1], [a_0;a_1,\dots,a_n] \big]$ or
$\big[ [a_0;a_1,\dots,a_n], [a_0;a_1,\dots,a_n+1] \big)$, depending on whether
$n$ is odd or even. Every half-open interval of this form obviously contains
an open interval, which proves the ``only if'' part of the lemma;
conversely, every open interval $(c,d)$ in $C$ contains a half-open interval of
this form (let $n-1$ be the first index at which the elements $c$ and $d$
differ), which proves the ``if'' part of the lemma.
\end{proof}

\begin{prop}
Let $A$ and $B$ be two countable dense subsets of $(0,1)$. Then there
exists an increasing homeomorphism $g:(0,1)\to(0,1)$ such that $g(A)=B$.
\label{gA.equals.B.prop}
\end{prop}

\begin{proof}
Note that any dense subset of $(0,1)$ must in fact be infinite. Fix
any well-orderings of $A$ and $B$ (that is, arrange the elements of $A$ and $B$
into infinite sequences); we emphasize that this well-ordering is not related
to the orderings of $A$ and $B$ as subsets of $(0,1)$. We recursively construct
a sequence of order-preserving bijections $g_j:A_j\to B_j$, where $A_j$ and
$B_j$ are subsets of $A$ and $B$, respectively, as follows. Choose any elements
$a_1\in A$ and $b_1\in B$, and define $A_1 = \{a_1\}$, $B_1=\{b_1\}$, and
$g_1(a_1)=b_1$.

If $n\ge2$ is even, we extend $A_{n-1}$, $B_{n-1}$, and $g_{n-1}$ as follows.
Choose the first (in the fixed well-ordering) element $b\in B\setminus
B_{n-1}$. If $b$ is smaller than every element of $B_{n-1}$, choose
$a\in A\setminus A_{n-1}$ that is smaller than every element of $A_{n-1}$. If
$b$ is larger than every element of $B_{n-1}$, choose
$a\in A\setminus A_{n-1}$ that is larger than every element of $A_{n-1}$.
If neither of these cases holds, then there are unique elements $c,d$ of
the finite set $B_{n-1}$ such that $c<b<d$ and $(c,d)\cap B_{n-1}=\emptyset$;
choose $a\in A\setminus A_{n-1}$ such that $g_{n-1}^{-1}(c) < a <
g_{n-1}^{-1}(d)$. (All of these choices are possible since $A$ is dense in
$(0,1)$.) After making this choice, we set $A_n = A_{n-1} \cup
\{a\}$ and $B_n = B_{n-1} \cup \{b\}$, and we define $g_n:A_n\to B_n$ by
$g_n(a)=b$ and $g_n(x)=g_{n-1}(x)$ if $x\in A_{n-1}$.

Similarly, if $n\ge3$ is odd, we extend $A_{n-1}$, $B_{n-1}$, and $g_{n-1}$ as
follows. Choose the first (in the fixed well-ordering) element $a\in A\setminus
A_{n-1}$. If $a$ is smaller than every element of $A_{n-1}$, choose
$b\in B\setminus B_{n-1}$ that is smaller than every element of $B_{n-1}$. If
$a$ is larger than every element of $A_{n-1}$, choose
$b\in B\setminus B_{n-1}$ that is larger than every element of $B_{n-1}$. If
neither of these cases holds, then there are unique elements $c,d$ of the
finite set $A_{n-1}$ such that $c<a<d$ and $(c,d)\cap A_{n-1}=\emptyset$;
choose $b\in B\setminus B_{n-1}$ such that $g_{n-1}(c) < b <
g_{n-1}(d)$. (All of these choices are possible since $B$ is dense in
$(0,1)$.) After making this choice, we set $A_n = A_{n-1} \cup \{a\}$
and $B_n = B_{n-1} \cup \{b\}$, and we define $g_n:A_n\to B_n$ by $g_n(a)=b$
and $g_n(x)=g_{n-1}(x)$ if $x\in A_{n-1}$.

It is easy to verify inductively that each $g_n$ is a bijection from
$A_n$ to $B_n$ that is order preserving with respect to the usual order on
$(0,1)$, and that $g_n|_{A_m} = g_m$ for all positive integers $m<n$.
Furthermore, the use of the well-orderings of $A$ and $B$ during the
construction forces $\bigcup_{n\ge1} A_n = A$ and $\bigcup_{n\ge1} B_n = B$.
Therefore, there is a unique function $g_\infty:A\to B$ (namely the union
of all the functions $g_n$) such that $g_\infty|_{A_n} = g_n$ for all $n\ge1$,
and in fact $g_\infty$ is an order-preserving bijection from $A$ to $B$.

Finally, define $g:(0,1)\to(0,1)$ by $g(x) = \sup
\{g_\infty(a)\colon a\in A,\, a<x\}$ for $x\in(0,1)$. Note that for any
$x\in(0,1)$, there exist $a_1,a_2\in A$ such that $a_1<x<a_2$ (by the
density of $A$), whence the set $\{g_\infty(a)\colon a\in A,\, a<x\}$ is
bounded above by
$g_\infty(a_1)<1$ and contains $g_\infty(a_2)>0$ by the order-preservation of
$g_\infty$; therefore $g(x)$ is a well-defined real number in $(0,1)$. Also,
$g$ is order preserving: if $c,d\in(0,1)$ with $c<d$, then there exists
$a_1,a_2\in A$ with $c<a_1<a_2<d$, whence
\begin{align*}
g(c) = \sup\{g_\infty(a)\colon a\in A,\, a<c\} &\le g_\infty(a_1) \\
&< g_\infty(a_2) \le \sup\{g_\infty(a)\colon a\in A,\, a<d\} = g(d).
\end{align*}
Moreover, $g$ is surjective: given $y\in(0,1)$, define
$x=\sup\{g_\infty^{-1}(b) \colon b\in B,\, b<y\}$. Because $g_\infty$ is an
order-preserving bijection from $A$ to $B$ and $B$ is dense near $y$, it is
easy to check that $g(x)=y$. Therefore $g$ is an order-preserving surjection
from $(0,1)$ to $(0,1)$, hence a homeomorphism by Lemma~\ref{orderbi.lemma}.
\end{proof}

A straightforward extension of this construction yields the following:

\begin{prop}
Let $A_1,A_2,\dots$ be a collection (finite or countably infinite) of pairwise
disjoint countable dense subsets of $(0,1)$, and similarly for $B_1,B_2,\dots$.
Then there exists an increasing homeomorphism $g:(0,1)\to(0,1)$ such that
$g(A_1)=B_1$, $g(A_2)=B_2$, and so on.
\label{countable.prop}
\end{prop}

\section{Proofs of the Theorems}

We call a function $h:\R\to\R$ a {\it chorus-line\/} function if $h$ maps
the set $(0,1)$ into itself and $h(x) = \floor x + h(\{x\})$ for all real
numbers $x$. This definition implies that $\floor{h(x)} = \floor x$ and
$\{h(x)\} = h(\{x\})$ for all $x$. Let $\G$ denote the set of increasing
homeomorphisms $g:(0,1)\to(0,1)$. To any $g\in\G$ we may associate a function
$\bar g:\R\to\R$, called the {\it chorus-line extension\/} of $g$, defined by
\begin{equation*}
\bar g(x) = \begin{cases}
x, &\text{if $x$ is an integer}, \\
\floor x + g(\{x\}), &\text{if $x$ is not an integer}.
\end{cases}
\end{equation*}
It is easy to see that $\bar g$ is an increasing homeomorphism from $\R$ to
itself and that both $\bar g$ and $\bar g^{-1} = \overline{g^{-1}}$ are
chorus-line functions. For any $g\in\G$, we define $f_g$ to be the restriction
of the function $g^{-1}\circ r\circ\bar g$ to the domain $(1,\infty)$, so that
$f_g(x) = g^{-1}(1/\bar g(x))$ for $x>1$. It is again easy to see that
$f_g\in\F$ with inverse $\phi_g = \bar g^{-1}\circ r\circ g$.

\begin{prop}
Let $g\in G$. Then $f_g$ is a faithful function satisfying
$\bar g([a_0;a_1,a_2,\dots]_{f_g}) = [a_0;a_1,a_2,\dots]_r$ for every
$[a_0;a_1,a_2,\dots]\in\C$. In other words, the continued $f$-expansion of
every real number $x$ is identical to the usual continued fraction expansion
of $\bar g(x)$.
\label{fg.faithful.prop}
\end{prop}

\begin{proof}
Let $x\in\R$, and let $E_r(\bar g(x))=[a_0;a_1,a_2,\dots]$; we want to show as
a first step that $E_{f_g}(x)=[a_0;a_1,a_2,\dots]$ as well. For each $n\ge0$,
define $y_n = [a_n;a_{n+1},a_{n+2},\dots]_r$ and $x_n = \bar g^{-1}(y_n)$, so
that $y_0=\bar g(x)$ and $x_0=x$. Notice that
\begin{equation}
\floor{x_n} = \floor{\bar g^{-1}(y_n)} = \floor{y_n} = a_n
\label{floor}
\end{equation}
for all $n\ge0$, where we have used the fact that $\bar g^{-1}$ is a
chorus-line function. Similarly,
\begin{align}
\phi_g(\{x_n\}) &= \bar g^{-1}\circ r\circ g (\{\bar g^{-1}(y_n)\}) \notag \\
&= \bar g^{-1}\circ r\circ g (\bar g^{-1}(\{y_n\})) \label{fracpart} \\
&= \bar g^{-1}\circ r(\{y_n\}) = \bar g^{-1}(y_{n+1}) = x_{n+1}. \notag
\end{align}
Now, simply considering the definition of the evaluation function
$E_{f_g}(x_0)$ in light of equations (\ref{floor}) and (\ref{fracpart})
reveals that $E_{f_g}(x) = E_{f_g}(x_0) = [a_0;a_1,a_2,\dots]$ as desired.

This shows that $E_{f_g} = E_r\circ\bar g$; in particular, $E_{f_g}$ is the
composition of two homeomorphisms and is therefore itself a homeomorphism, with
$E_{f_g}^{-1} = \bar g^{-1}\circ E_r^{-1}$. To show that $f_g$ is faithful, it
only remains to show that $V_{f_g}$ is well-defined on all of $\C$ and inverts
the function $E_{f_g}$. As above, let $x\in\R$ and $E_{f_g}(x)=E_r(\bar
g(x))=[a_0;a_1,a_2,\dots]$; we want to show that $V_{f_g}([a_0;a_1,a_2,\dots])
= \lim_{n\to\infty} V_{f_g}([a_0;a_1,\dots,a_n])$ exists and equals $x$. Now
$V_{f_g}$ does invert $E_{f_g}$ on $\C_t$, so
\begin{align*}
V_{f_g}([a_0;a_1,\dots,a_n]) &= E_{f_g}^{-1}([a_0;a_1,\dots,a_n]) \\
&= \bar g^{-1}\circ E_r^{-1}([a_0;a_1,\dots,a_n]) \\
&= \bar g^{-1}\circ V_r([a_0;a_1,\dots,a_n]) = \bar g^{-1}
([a_0;a_1,\dots,a_n]_r).
\end{align*}
However, we know that $\lim_{n\to\infty} [a_0;a_1,\dots,a_n]_r =
[a_0;a_1,a_2,\dots]_r$ by the convergence of usual continued fractions.
Therefore, by the continuity of $\bar g^{-1}$, we have
\begin{equation*}
V_{f_g}([a_0;a_1,a_2,\dots]) = \lim_{n\to\infty} \bar g^{-1}
([a_0;a_1,\dots,a_n]_r) = \bar g^{-1}([a_0;a_1,a_2,\dots]_r) = \bar g^{-1}(\bar
g(x)) = x
\end{equation*}
as desired.
\end{proof}

Together, Propositions \ref{gA.equals.B.prop}--\ref{fg.faithful.prop} imply
each of the Theorems 1--4. For example, given $x,y\in(0,1)$, let $g\in\G$ be
chosen so that $g(x)=y$. Then Proposition \ref{fg.faithful.prop} tells us that
$f_g$ is a faithful function in $\F$ and that $E_{f_g}(x) = E_r(y)$, which is
precisely the statement of Theorem 1. In fact, $g$ can be chosen to be
piecewise linear, in which case $f_g$ is given piecewise by M\"obius
transformations $(ax+b)/(cx+d)$.

For any positive integer $d$, let $Q(d)$ denote the set of numbers in $(0,1)$
that are algebraic over $\Q$ of degree exactly $d$, so that $Q(1)=\Q\cap(0,1)$
for instance. Each $Q(d)$ is a countable dense subset of $(0,1)$. Therefore,
Proposition \ref{gA.equals.B.prop} tells us that there exists a function
$g\in\G$ such that $g(Q(3))=Q(2)$. We then know from Proposition
\ref{fg.faithful.prop} that the corresponding $f_g$ has the property that the
continued $f_g$-expansion of any number in $Q(3)$ is the same as the usual
continued fraction expansion of a number in $Q(2)$, and vice versa. Since the
usual continued fraction expansion of a real number is periodic if and only if
the number is a quadratic irrational, this $f_g$ gives a faithful function
such that the continued $f_g$-expansion of a real number $x$ is periodic if and
only if $x$ is a cubic irrational number, establishing Theorem~\ref{cubic.thm}.

A similar approach establishes Theorem \ref{frage.thm}. In fact, the
singular Minkowski function $?(x)$ (see \cite{frage}) is an increasing
homeomorphism of $(0,1)$ that was constructed to have the property that
$?(Q(1)\cup Q(2))=Q(1)$. Therefore $f_?\in\F$ is a faithful function such that
the continued $f_?$-expansion of any rational or quadratic irrational number
terminates.

As for Theorem \ref{simultaneous.thm}, we partition $\C_t$ into infinitely
many sets $\C_t(1),\C_t(2),\dots$, where we define
$$
\C_t(d) = \{ [0; a_1, \dots, a_n] \colon n\ge1,\, \text{each }a_j\in\Z,\,
a_j\ge1 \text{ for all } 1\le j<n,\, a_n=d+1 \}.
$$
Each $\C_t(d)$ is dense in $\C$ by Lemma \ref{dense.lemma}. Therefore, if we
define $Q(1,d) = V_r(\C_t(d))$, then $Q(1,1),Q(1,2),\dots$ is a partition of
$\Q\cap(0,1)$ into countably many countably infinite subsets, each dense in
$(0,1)$. Applying Proposition \ref{countable.prop}, we can find a function
$g\in\G$ such that $g(Q(d))=Q(1,d)$ for every positive integer $d$. Then, by
Proposition \ref{fg.faithful.prop}, the set of real numbers in $(0,1)$ whose
continued $f_g$-expansion terminates in the integer $d+1$ is precisely $Q(d)$.

We briefly discuss the extension of Theorem \ref{simultaneous.thm}
mentioned in the introduction. We generalize the notation of the previous
paragraph by defining $\C_t(n_1,n_2,\dots,n_k)$ to be the set of
terminating expansions in $\C$ that begin with a zero and end with the $k$
integers $n_1,\dots,n_k$, and we set $Q(1,n_1,\dots,n_k) =
V_r(\C_t(n_1,n_2,\dots,n_k))$. We choose a function that encodes every
integer, positive or negative, as a positive integer; one such function is
\begin{equation*}
\eta(k) = \begin{cases}
2|k|+1, &\text{if }k\le0; \\ 2k, &\text{if }k>0,
\end{cases}
\end{equation*}
which is a bijection from $\Z$ to $\Z^+$ whose inverse is $\delta(k) = (-1)^k
\floor{\frac k2}$. Using a modification of Proposition \ref{countable.prop},
we can find a function $g\in\G$ such that, for each algebraic number
$x\in(0,1)$, if the minimal polynomial of $x$ is $c_dt^d + c_{d-1}t^{d-1} +
\dots + c_1t + c_0$, then
$g(x)\in\C_t(\eta(c_0),\eta(c_1),\dots,\eta(c_{d-1}),\eta(c_d),d+1)$. (Many
transcendental numbers would also be mapped into these subsets of $\C_t$ by
$g$.) The corresponding faithful function $f_g\in\F$ would then have the
property that the continued $f_g$-expansion of $x$ would terminate in the
sequence $\eta(c_0),\eta(c_1),\dots,\eta(c_{d-1}),\eta(c_d),d+1$, thus
encoding the minimal polynomial of $x$. If we had an oracle that could compute
this function $f_g$ quickly, we could test whether any real number $y$ was
algebraic by computing its continued $f_g$-expansion; if it terminated, say as
$[0;a_1,\dots,a_n]$, then $y$ would be either transcendental or else a root of
the polynomial
\begin{equation*}
\delta(a_{n-1})t^{a_n-1} + \delta(a_{n-2})t^{a_n-2} + \dots +
\delta(a_{n-(a_n-1)})t + \delta(a_{n-a_n}).
\end{equation*}
Of course, this is only a fantasy, as the function $f_g$ is hopelessly
infeasible for exact computation. Other types of encoding functions are
possible, of course; for example, one can encode every finite sequence of
integers as a single positive integer via some G\"odel-type code.

Proposition \ref{fg.faithful.prop} shows us how we can contruct a faithful
function in $\F$ from any function in $\G$. The following result demonstrates
that the opposite is also true:

\begin{prop}
Let $f\in\F$ be a faithful function. Define $g$ to be the restriction of the
function $V_r\circ E_f$ to $(0,1)$. Then $g\in\G$ and $f_g=f$.
\label{cl.conjugate.prop}
\end{prop}

\begin{proof}
Since $E_f:(0,1)\to\C^0$ and $V_r:\C^0\to(0,1)$ are both order preserving,
we see from Lemma \ref{orderbi.lemma} that $g$ is indeed a homeomorphism from
$(0,1)$ to itself, and hence $g\in\G$. Notice that $V_r\circ E_f$ has the
property that $V_r\circ E_f(x) = \floor x + V_r\circ E_f(\{x\})$; therefore
$V_r\circ E_f$ is its own chorus-line extension and so $\bar g=V_r\circ E_f$.
If we write $E_f(x) = [a_0; a_1, a_2, \dots]$, then $\bar g(x) = V_r([a_0;
a_1, a_2, \dots]) = [a_0; a_1, a_2, \dots]_r$ and therefore $r\circ\bar g(x) =
[0; a_0, a_1, a_2,
\dots]_r$. Finally, note that $(V_r\circ E_f)^{-1} = E_f^{-1}\circ V_r^{-1} =
V_f\circ E_r$, and so
\begin{align*}
g^{-1}\circ r\circ\bar g(x) = V_f\circ E_r([0; a_0, a_1, a_2, \dots]_r) &= [0;
a_0, a_1, a_2, \dots]_f \\
&= 0+f([a_0; a_1, a_2, \dots]_f) = f(x)
\end{align*}
as desired.
\end{proof}

We say that two functions $f,f'\in\F$ are {\it chorus-line conjugate\/} if
there exists an increasing continuous chorus-line function $h$ such that $f$
equals the restriction of $h^{-1}\circ f'\circ h$ to the domain $(1,\infty)$.
It is easy to check that chorus-line conjugacy is an equivalence relation on
$\F$. The next theorem shows that the equivalence class containing the
reciprocal function $r$ is precisely the class of faithful functions.

\begin{theorem}
Let $f\in\F$. The following are equivalent:
\begin{enumerate}
\renewcommand{\theenumi}{{\it\roman{enumi}}}
\item $f$ is faithful;
\item $f$ is chorus-line conjugate to $r$;
\item $f=f_g$ for some $g\in\G$.
\end{enumerate}
In particular, there is a one-to-one correspondence between functions
$g\in\G$ and faithful functions $f\in\F$.
\label{fisg.thm}
\end{theorem}

\begin{proof}
It is easy to see that any increasing continuous chorus-line function is equal
to $\bar g$ for some $g\in\G$. Thus $f$ and $r$ are chorus-line conjugate
if and only if $f=\bar g^{-1}\circ r\circ\bar g=f_g$ for some $g\in\G$, which
shows the equivalence of statements (ii) and (iii). Proposition
\ref{fg.faithful.prop} shows that statement (iii) implies statement (i), while
Proposition \ref{cl.conjugate.prop} shows that statement (i) implies
statement (iii). Therefore the three statements are indeed equivalent.

The assertion that there is a one-to-one correspondence between functions
$g\in\G$ and faithful functions $f\in\F$ requires some justification, as it is
not immediately clear that different functions $g,g'\in\G$ give rise to
distinct $f_g$ and $f_{g'}$. Suppose that $g\ne g'$, and choose an $x\in(0,1)$
such that $g(x)\ne g'(x)$. By Proposition \ref{fg.faithful.prop}, we have
$E_{f_g}(x) = E_r(g(x)) \ne E_r(g'(x)) = E_{f_{g'}}(x)$. Since the continued
$f$- and $f'$-expansions of $x$ differ, we must have $f_g\ne f_{g'}$.
\end{proof}

At the end of the last proof, we used the fact that two faithful functions
whose corresponding expansion functions are different must be distinct. In
fact, the converse is also true:

\begin{prop}
Suppose that $f_1,f_2\in\F$ are faithful functions with the property that
$E_{f_1}=E_{f_2}$. Then $f_1=f_2$.
\label{converse.prop}
\end{prop}

\begin{proof}
Given $x\in(1,\infty)$, we want to prove that $f_1(x)=f_2(x)$. Write
$E_{f_1}(x) = [a_0; a_1, a_2, \dots] = E_{f_2}(x)$. Note that the continued
$f_1$-expansion of $f_1(x) = 0+f_1(x)$ is $[0; a_0, a_1, a_2, \dots]$, so that
$E_{f_1}(f_1(x)) = [0; a_0, a_1, a_2, \dots]$. Similarly, $E_{f_2}(f_2(x)) =
[0; a_0, a_1, a_2, \dots]$. But then $E_{f_1}(f_2(x)) = [0; a_0, a_1, a_2,
\dots]$ since $E_{f_1}=E_{f_2}$. Therefore $E_{f_1}(f_1(x)) =
E_{f_1}(f_2(x))$, and since $f_1$ is faithful, $E_{f_1}$ is injective and
thus $f_1(x)=f_2(x)$.
\end{proof}

\section{Continued Power Function Expansions}

Let us consider a particular one-parameter family of functions from $\F$,
namely the power functions $f_\alpha(x) = x^{-\alpha}$ for $\alpha>0$, so that
$f_1=r$. A continued $f_\alpha$-expansion of a real number is thus an
expression of the form
\begin{equation*}
a_0 + (a_1 + (a_2 + (a_3 + (a_4 + \cdots
)^{-\alpha})^{-\alpha})^{-\alpha})^{-\alpha}.
\end{equation*}
At a problem session of the West Coast Number Theory Conference in 1999, Kevin
O'Bryant considered the case $\alpha=1/2$, which he called the continued root
expansion of a real number~$x$. For instance, some rational numbers such as
\begin{equation*}
\frac23 = 0 + \frac1{\sqrt{2 + \frac1{\sqrt{16}}}} = [0;2,16]_{f_{1/2}}
\end{equation*}
and $\frac{27}{47} = [0;3,1098,2892,410,256]_{f_{1/2}}$ have terminating
continued root expansions. On the other hand, O'Bryant remarked that
\begin{equation*}
E_{f_{1/2}}(\tfrac34) = [0; 1, 1, 2, 8, 5, 1, 3, 3, 14, 321, 2, 300, 1, 13, 2,
6, 1, 1, 2,\dots]
\end{equation*}
does not seem to terminate; but we do not know how to prove this. At the same
problem session, Bart Goddard noted several other examples; for instance,
$E_{f_5}(\root5\of7) = [1;1,1,1,\dots]$, and the continued $f_{3/2}$-expansion
of $\root3\of3 = 1.44224957$ looks at first to be periodic of period four.
However,
$$
E_{f_{3/2}}(\root3\of3) =
[1;1,1,2,1,1,1,2,1,1,1,2,1,1,1,3,1,1,1,1,3,1,2,1,1,7,23,1,\dots]
$$
does not seem to be periodic, while
$$
x = V_{f_{3/2}}([1;1,1,2,1,1,1,2,1,1,1,2,\dots]) = 1.44225029
$$
is the nearby number that satisfies the equation $x = 1 + (1 + (1 + (2 +
x^{-3/2} )^{-3/2} )^{-3/2})^{-3/2}$. In fact, writing $w$, $z$, and $y$ as the
three quantities in parentheses on the right-hand side of this equation
(starting from the innermost parentheses), we see that this number $x$ is the
$x$-coordinate of one solution to the system of equations
\begin{equation*}
y^3(x-1)^2=1, \quad z^3(y-1)^2=1, \quad w^3(z-1)^2=1, \quad x^3(w-2)^2=1.
\end{equation*}
Using elimination theory and a computational algebra package, we can show that
this number $x$ is algebraic of degree 93; more precisely, it is the fourth of
seven real roots of the irreducible polynomial

{\medskip\tiny\noindent\tolerance=100 \emergencystretch=60pt
${-2401} + 12348x - 22442x^2 + 275800x^3 - 1337555x^4 +
2423872x^5 - 15418480x^6 + 70540444x^7 - 127417629x^8 +
557491285x^9 - 2405709582x^{10} + 4329064154x^{11} -
14625356403x^{12} + 59525595995x^{13} - 106704972668x^{14} +
296336967716x^{15} - 1137325584809x^{16} + 2031978559593x^{17} -
4823156208926x^{18} + 17439240838410x^{19} - 31080157671439x^{20} +
64755935263191x^{21} - 220128009411364x^{22} + 391629168869836x^{23} -
730457802870121x^{24} + 2326819690217101x^{25} -
4133151272936538x^{26} + 7008134858873413x^{27} -
20830642337065947x^{28} + 36915446021983793x^{29} -
57610938763130172x^{30} + 159042225095378801x^{31} -
280663409776128761x^{32} + 407135007678293093x^{33} -
1039187243118822998x^{34} + 1820679151953897429x^{35} -
2473096725871085456x^{36} + 5813652442623811749x^{37} -
10072679129493093706x^{38} + 12871647850762706121x^{39} -
27781471420314300292x^{40} + 47381905470113399929x^{41} -
57050802593263213282x^{42} + 112791554026161912586x^{43} -
188441383925133877380x^{44} + 213305384930045048629x^{45} -
385514084983483335018x^{46} + 627994700913994360904x^{47} -
663244567700776798728x^{48} + 1093132419527821059119x^{49} -
1729925336326989733586x^{50} + 1677221045253645425438x^{51} -
2509661817013886167565x^{52} + 3855357757564686581031x^{53} -
3317198776437539175677x^{54} + 4459703033666698711524x^{55} -
6699690356527027205257x^{56} + 4712835619541841129045x^{57} -
5515895022295988234329x^{58} + 8418986026761406210570x^{59} -
3557493739786303419709x^{60} + 2967144052401077482669x^{61} -
6097230519934320899607x^{62} - 2428850299115623275704x^{63} +
4729361127147513235131x^{64} - 963376378662124156885x^{65} +
11586996393861391188069x^{66} - 14002890770158170207742x^{67} +
7817852618475056747791x^{68} - 16575420482318059675255x^{69} +
16734537138028434957244x^{70} - 7089221324585019485886x^{71} +
10933203152887989317941x^{72} - 8907308593248131984589x^{73} -
1468466124295786901110x^{74} + 1320989510532080943648x^{75} -
2515839920964633664993x^{76} + 8733255119045045834197x^{77} -
8648102018082906320368x^{78} + 7228474474951901475700x^{79} -
7991643882573751006683x^{80} + 6679169691105510026567x^{81} -
4448501164546530714930x^{82} + 3131234600047636654702x^{83} -
1971348622249197779737x^{84} + 953049356660824435629x^{85} -
363152194705059550764x^{86} + 72557323790158601616x^{87} +
29118868029709313904x^{88} - 24792096645669431805x^{89} +
4762766285696524504x^{90} + 768224935099977754x^{91} -
343187952655081548x^{92} + 28598996054590129x^{93}$.

\medskip}

Continued $f_\alpha$ expansions exhibit several interesting phenomena which
merit further study. For example, let us consider whether the limit defining
$[1;1,1,1,\dots]_{f_\alpha}$ converges. Computationally, we find that there
is a threshhold number $\alpha_0 = 4.1410415\dots$ with the property that
$[1;1,1,1,\dots]_{f_\alpha}$ converges for all $0<\alpha<\alpha_0$ (which
we refer to as ``small $\alpha$'') but diverges for all $\alpha>\alpha_0$
(``large $\alpha$''). In fact, $\alpha_0$ is the unique positive solution of
the equation $y^y=(y-1)^{y+1}$.

The behavior of $E_{f_\alpha}$ as $\alpha$ passes through $\alpha_0$
experiences a classic bifurcation. For all positive $\alpha$, the function
$(x-1)^{-1/\alpha}$ has a unique fixed point between 1 and 2; however, this
fixed point is repelling for small $\alpha$ but attracting for large $\alpha$.
Therefore only a single real number has $[1;1,1,1,\dots]$ as its continued
$f_\alpha$-expansion for small $\alpha$, and $V_{f_\alpha}[1;1,1,1,\dots]$
converges back to this real number. In contrast, there is a whole interval of
real numbers having $[1;1,1,1,\dots]$ as their continued $f_\alpha$-expansions
for large $\alpha$. For example, when $\alpha=5$, all real numbers in the
interval $(1.06377, 1.73411)$ have $[1;1,1,1,\dots]$ as their continued
$f_5$-expansion. In particular, since $\root5\of7=1.47577$, the example
$E_{f_5}(\root5\of7) = [1;1,1,1,\dots]$ mentioned above is less significant
than it seems.

Indeed, for large $\alpha$ the function $((x-1)^{-\alpha}-1)^{-\alpha}$ has
three fixed points in the interval $(1,2)$, the central one being unstable and
the outer two being stable. The evaluations $V_{f_\alpha}[1;1,1,\dots,1]$
oscillate back and forth between ever-decreasing neighborhoods of the two
stable fixed points as the number of ones increases, and hence
$V_{f_\alpha}[1;1,1,1,\dots]$ does not converge for large $\alpha$. The two
outer fixed points approach 1 from above and 2 from below, respectively, as
$\alpha$ tends to infinity. We can conclude, for instance, that for every real
number $x\in(1,2)$, there exists an $\alpha(x)$ such that, whenever
$\alpha>\alpha(x)$, we have $E_{f_\alpha}(x) = [1;1,1,1,\dots]$. For example,
we have $E_{f_\alpha}(\root5\of7) = [1;1,1,1,\dots]$ for all $\alpha>4.26159$.

The above discussion implies in particular that $f_\alpha$ is not faithful for
large $\alpha$. On the other hand, it can be shown that no analogous
bifurcation occurs for $[n;n,n,n,\dots]$ when $n\ge2$ (the key equality now
becomes $n^{y+1}y^y=(y-1)^{y+1}$ which has no positive solution).
Computational evidence suggests that periodic sequences of longer period never
undergo bifurcations either. We are thus led to conjecture that $f_\alpha$ is
faithful for all $0<\alpha<\alpha_0$. (Kakeya's theorem (see \cite[Section
8.3]{schweiger}) is only relevant when $\alpha<1$.) In particular, this
conjecture would imply by Theorem \ref{fisg.thm} that the power
function $x^{-\alpha}$ is chorus-line-conjugate to the reciprocal function on
$(1,\infty)$ for $\alpha<\alpha_0$ but not for $\alpha>\alpha_0$, a curious
state of affairs.

We can also use these functions to show that the analogue of Proposition
\ref{converse.prop} for non-faithful functions does not hold. Indeed,
if $f\in F$ is any function agreeing with $f_5$ outside the interval
$(1.06377,1.73411)$, then it is easy to see that $E_f = E_{f_5}$.

We mention one last phenomenon, where we fix a real number $x$ and consider
the function from $\R$ to $\C$ that maps $\alpha$ to $E_{f_\alpha}(x)$.
Counterintuitively, this function is not an order-preserving function of
$\alpha$. For example, when $x=\frac12$, the function $\alpha\mapsto
E_{f_\alpha}(\frac12)$ is increasing for $\alpha<2.24228$ but decreasing
thereafter (stabilizing eventually at $[1;1,1,1,\dots]$, as we have already
seen).

\bigskip
{\noindent\small{\it Acknowledgements.} The author acknowledges the support
of the Natural Sciences and Engineering Research Council and of the
Department of Mathematics of the University of British Columbia, and thanks
Bart Goddard for the interesting questions.}

\bibliographystyle{amsplain}
\bibliography{unreasonable}

\end{document}